\documentclass[11pt]{article}

\usepackage{amsmath,amssymb}

\makeatletter
\@addtoreset{equation}{section}
\makeatother

\title{\bf q-ANALOGUE OF THE BEREZIN QUANTIZATION METHOD}

\author{\sl D. Shklyarov \and S. Sinel'shchikov \and \sl L. Vaksman}

\date{\tt Institute for Low Temperature Physics \& Engineering \\
National Academy of Sciences of Ukraine}

\newtheorem{theorem}{Theorem}[section]
\newtheorem{lemma}[theorem]{Lemma}
\newtheorem{proposition}[theorem]{Proposition}

\begin{document}
\large
\maketitle

\bigskip

\section{Introduction}

 The concept of deformation quantization of a symmetric manifold $M$ has
been defined by Bayen, Flato, Fronsdal, Lichnerovich, and Sternheimer in
\cite{BFFLS}. Deformation quantization means a formal $*$-product
$$f_1*f_2=f_1 \cdot f_2+\sum_{k=1}^\infty C_k(f_1,f_2)t^k,\qquad f_1,f_2 \in
C^\infty(M)$$
with some additional properties, where $C_k:C^\infty(M)\times C^\infty(M)\to
C^\infty(M)$ are bidifferential operators.

  In the special case of the unit disc in ${\mathbb C}$ with
$SU_{1,1}$-invariant symplectic structure, a formal $*$-product and explicit
fopmulae for $C_k$, $k \in{\mathbb N}$, are derivable by a method of Berezin
\cite{CGR, B}.

 Our intention is to replace the ordinary disc with its q-analogue. We
are going to produce $U_q \mathfrak{su}_{11}$-invariant formal deformation
of our quantum disc and to obtain an explicit formula for $C_k$, $k
\in{\mathbb N}$, using a q-analogue of the Berezin method \cite{SSV5}.

 Our work is closely related to the paper of Klimek and Lesniewski \cite{KL}
on two-parameter deformation of the unit disc. The explicit formulae for
$C_k$ we provide below work as a natural complement to the results of this
paper.

\bigskip

\section{Covariant symbols of linear operators}

 Everywhere in the sequel the field of complex numbers ${\mathbb C}$ is
assumed as a ground field. Let also $q \in(0,1)$.

 Consider the well known algebra ${\rm Pol}({\mathbb C})_q$ with two
generators $z,z^*$ and a single commutation relation $z^*z=q^2zz^*+1-q^2$.
Our intention is to produce a formal $*$-product
\begin{equation}f_1*f_2=f_1 \cdot f_2+\sum_{k=1}^\infty
C_k(f_1,f_2)t^k,\qquad f_1,f_2 \in{\rm Pol}({\mathbb C})_q,
\end{equation}
(with some remarkable properties) to be given by explicit formulae for
bilinear operators $C_k:{\rm Pol}({\mathbb C})_q \times{\rm Pol}({\mathbb
C})_q \to{\rm Pol}({\mathbb C})_q$.

 We describe in this section the method of producing this $*$-product whose
idea is due to F. Berezin.

 It was explained in \cite{SSV1} that the vector space $D({\mathbb U})_q'$
of formal series $\displaystyle \sum_{j,k=0}^\infty a_{jk}z^jz^{*k}$ with
complex coefficients is a q-analogue of the space of distributions in the
unit disc ${\mathbb U}=\{z \in{\mathbb C}|\;|z|<1 \}$. Equip this space of
formal series with the topology of coefficientwise convergence. Since
$\{z^jz^{*k}\}_{j,k \in{\mathbb Z}_+}$ constitute a basis in the vector
space ${\rm Pol}({\mathbb C})_q$, ${\rm Pol}({\mathbb C})_q$ admits an
embedding into $D({\mathbb U})_q'$ as a dense linear subvariety.

 Consider the unital subalgebra ${\mathbb C}[z]_q \subset{\rm Pol}({\mathbb
C})_q$ generated by $z \in{\rm Pol}({\mathbb C})_q$. Let $\alpha>0$. We
follow \cite{KL} in equipping the vector space ${\mathbb C}[z]_q$ with the
scalar product
$(z^j,z^k)_\alpha=\delta_{jk}\dfrac{(q^2;q^2)_k}{(q^{4 \alpha+2};q^2)_k}$,
$j,k \in{\mathbb Z}_+$, where
$(a;q^2)_k=(1-a)(1-q^2a)\ldots(1-q^{2(k-1)}a)$. Let $L_a^2(d \nu_ \alpha)_q$
be the a completion of ${\mathbb C}[z]_q$ with respect to the norm $\| \psi
\|_\alpha=(\psi,\psi)_\alpha^{1/2}$. It was demonstrated in \cite{KL} that
the Hilbert space $L_a^2(d \nu_ \alpha)_q$ is a q-analogue of the weighted
Bergman space. Let $\widehat{z}$ be a linear operator of
multiplication by $z$:
$$\widehat{z}:L_a^2(d \nu_ \alpha)_q \to L_a^2(d \nu_ \alpha)_q;\qquad
\widehat{z}:\psi(z)\mapsto z \cdot \psi(z),$$
and denote by $\widehat{z}^*$ the adjoint operator in $L_a^2(d \nu_
\alpha)_q$ to $\widehat{z}$. The
definition of the scalar product in $L_a^2(d \nu_ \alpha)_q$ implies that
the operators $\widehat{z}$, $\widehat{z}^*$ are bounded.
 Equip the space ${\cal L}_\alpha$ of bounded
linear operators in $L_a^2(d \nu_ \alpha)_q$ with the weakest topology in
which all the linear functionals
$$l_{\psi_1,\psi_2}:{\cal L}_\alpha \to{\mathbb C},\qquad
l_{\psi_1,\psi_2}:A \mapsto(A \psi_1,\psi_2)_\alpha,\qquad \psi_1,\psi_2
\in{\mathbb C}[z]_q$$
are continuous. The following proposition is a straightforward consequence
of the definitions (see the proof in \cite{SSV1}).

\medskip

\begin{proposition} Given any bounded linear operator $\widehat{f}$ in the
Hilbert space $L_a^2(d \nu_ \alpha)_q$, there exists a unique formal series
$f=\sum \limits_{j,k=0}^\infty a_{jk}z^jz^{*k}\in D({\mathbb U})_q'$ such
that $\widehat{f}=\sum \limits_{j,k=0}^\infty a_{jk}\widehat{z}^j
\widehat{z}^{*k}$.
\end{proposition}

\medskip

 Thus we get an injective linear map ${\cal L}_\alpha \to D({\mathbb
U})_q'$, $\widehat{f}\mapsto f$. The distribution $f$ is called a {\sl
covariant symbol} of the linear operator $\widehat{f}$.

\medskip

 {\sc Remark 2.2.} For an arbitrary  $f \in{\rm Pol}({\mathbb
C})_q$, there exists a unique operator $\widehat{f}\in{\cal L}_\alpha$ with
the covariant symbol $f$. Specifically, for $f=\sum
\limits_{j,k=0}^{N(f)}a_{jk}z^jz^{*k}$, one has $\widehat{f}=\sum
\limits_{j,k=0}^{N(f)}a_{jk}\widehat{z}^j \widehat{z}^{*k}$.

\medskip

 We follow F. Berezin in producing the $*$-product of covariant symbols
using the ordinary product of the associated linear operators.

 Let $f_1,f_2 \in{\rm Pol}({\mathbb C})_q$ and $\widehat{f}_1,\widehat{f}_2
\in{\cal L}_\alpha$ be the operators whose covariant symbols are $f_1$,
$f_2$. Under the notation $t=q^{4 \alpha}$, let $m_t(f_1,f_2)$ stand for the
covariant symbol of the product $\widehat{f}_1 \cdot \widehat{f}_2$ of the
linear maps $\widehat{f}_1$, $\widehat{f}_2$. Evidently, we have constructed
a bilinear map $m_t:{\rm Pol}({\mathbb C})_q \times{\rm Pol}({\mathbb C})_q
\to D({\mathbb U})_q'$.

 The $*$-product $f_1*f_2$ of $f_1,f_2 \in{\rm Pol}({\mathbb C})_q$ is to be
introduced by replacement of the one-parameter family of distributions
$m_t(f_1,f_2)$, $t \in(0,1)$, with its asymptotic expansion as $t \to 0$.

\bigskip

\section{\boldmath$*$-Product}

 The term {\sl 'order one differential calculus over the algebra ${\rm
Pol}({\mathbb C})_q$'} stand for a ${\rm Pol}({\mathbb C})_q$-bimodule
$\Omega^1({\mathbb C})_q$ equipped with a linear map $d:{\rm Pol}({\mathbb
C})_q \to \Omega^1({\mathbb C})_q$ such that \\ i) $d$ satisfies the
Leibniz rule $d(f_1f_2)=df_1 \cdot f_2+f_1 \cdot df_2$ for any $f_1,f_2
\in{\rm Pol}({\mathbb C})_q$,\\ ii) $\Omega^1({\mathbb C})_q$ is a linear
span of $f_1 \cdot df_2 \cdot f_3$, $f_1,f_2,f_3 \in{\rm Pol}({\mathbb
C})_q$ (see \cite{KS}).

 One can find in \cite{SV} a construction of that kind of order one
differential calculus for a wide class of prehomogeneous vector spaces $V$.
In the case $V={\mathbb C}$ we deal with this calculus is well known; it can
be described in terms of the following commutation relations:
$$z \cdot dz=q^{-2}dz \cdot z,\qquad z^*dz^*=q^2dz^*z^*,\qquad z^*dz=q^2dz
\cdot z^*,\qquad z \cdot dz^*=q^{-2}dz^*z.$$

 The partial derivatives $\dfrac{\partial^{(r)}}{\partial z}$,
$\dfrac{\partial^{(r)}}{\partial z^*}$, $\dfrac{\partial^{(l)}}{\partial
z}$, $\dfrac{\partial^{(l)}}{\partial z^*}$ are
linear operators in ${\rm Pol}({\mathbb C})_q$ given by
$$df=\frac{\partial^{(r)}f}{\partial z}dz+\frac{\partial^{(r)}f} {\partial
z^*}dz^*=dz \frac{\partial^{(l)}f}{\partial
z}+dz^*\frac{\partial^{(l)}f}{\partial z^*},$$
with  $f \in{\rm Pol}({\mathbb C})_q$.

 Let $\widetilde{\square}:{\rm Pol}({\mathbb C})_q^{\otimes 2}\to{\rm
Pol}({\mathbb C})_q^{\otimes 2}$, $m_0:{\rm Pol}({\mathbb C})_q^{\otimes
2}\to{\rm Pol}({\mathbb C})_q$ be linear operators given by
$$\widetilde{\square}(f_1 \otimes
f_2)=\left(\frac{\partial^{(r)}f_1}{\partial z^*}\otimes 1 \right)\cdot
q^{-2}(1-(1+q^{-2})z^*\otimes z+q^{-2}z^{*2}\otimes z^2)\cdot \left(1
\otimes \frac{\partial^{(l)}f_2}{\partial z}\right),$$
$m_0(f_1 \otimes f_2)=f_1f_2$, with $f_1,f_2 \in{\rm Pol}({\mathbb C})_q$.

\medskip

\begin{theorem}\label{ae}For all $f_1,f_2 \in{\rm Pol}({\mathbb C})_q$, the
following asymptotic expansion in $D({\mathbb U})_q'$ is valid:
$$m_t(f_1,f_2)\sim_{_{\!\!\!\!\!\!\!\!\!t \to 0}}f_1*f_2,\qquad{\rm with}$$
\begin{equation}\label{stpr}f_1*f_2=f_1 \cdot f_2+\sum_{k=1}^\infty
C_k(f_1,f_2)t^k \in{\rm Pol}({\mathbb C})_q[[t]],
\end{equation}
\begin{equation}\label{Ck}C_k(f_1,f_2)=m_0 \left(\left(p_k
\left(\widetilde{\square}\right)-
p_{k-1}\left(\widetilde{\square}\right)\right)(f_1 \otimes f_2)\right),
\end{equation}
and $p_k(x)$, $k \in{\mathbb Z}_+$, are polynomials given by
\begin{equation}\label{pk}p_k(x)=\sum_{j=0}^k
\frac{(q^{-2k};q^2)_j}{(q^2;q^2)^2_j}q^{2j}
\prod_{i=0}^{j-1}(1-q^{2i}((1-q^2)^2 x+1+q^2)+q^{4i+2}).
\end{equation}
\end{theorem}

\medskip

 This statement is to be proved in the next section, using the results of
\cite{SSV5} on a q-analogue of the Berezin transform \cite{UU}.

 We are grateful to H. T. Koelink who attracted our attention to the fact
that the polynomials $p_k(x)$ differ from the polynomials of Al-Salam --
Chihara \cite{K} only by normalizing multiples and a linear change of the
variable $x$.

\bigskip

\section{A q-analogue of the Berezin transform}

 Remind the notation $t=q^{4 \alpha}$, with $q \in(0,1)$, $\alpha>0$.

 Consider the linear map ${\rm Pol}({\mathbb C})_q \to{\cal L}_\alpha$ which
sends a polynomial $\stackrel{\circ}{f}=\sum \limits_{jk}b_{jk}z^{*j}z^k$ to
the linear operator $\widehat{f}=\sum
\limits_{jk}b_{jk}\widehat{z}^{*j}\widehat{z}^k$. The polynomial
$\stackrel{\circ}{f}$ will be called a {\sl contravariant symbol} of the
linear operator $\widehat{f}$.

 Note that our definitions of covariant and contravariant symbols agree with
the conventional ones, as one can observe from \cite{SSV5} (specifically,
see proposition 6.6 and lemma 7.2 of that work).

 The term {\sl 'q-transform of Berezin'} will be stand for the linear
operator $B_{q,t}:{\rm Pol}({\mathbb C})_q \to D({\mathbb U})_q'$,
$B_{q,t}:\stackrel{\circ}{f}\mapsto f$, which sends the contravariant
symbols of linear operators $\widehat{f}=\sum
\limits_{jk}b_{jk}\widehat{z}^{*j}\widehat{z}^k$ to their covariant symbols.

\medskip

 {\sc Remark 4.1.} It is easy to extend the operators $B_{q,t}$ onto the
entire {\sl 'space of bounded functions in the quantum disc'} via a
non-standard approach to their construction (see \cite{SSV5}).

\medskip

\cite[proposition 5.5]{SSV5} imply

\medskip

\begin{proposition}\label{Bqtae}Given arbitrary $\stackrel{\circ}{f}\in{\rm
Pol}({\mathbb C})_q$, the following {\sl asymptotic expansion} in the
topological vector space $D({\mathbb U})_q'$ is valid:
$$B_{q,t}\stackrel{\circ}{f}\;\sim_{_{\!\!\!\!\!\!\!\!\!t \to
0}}\;\stackrel{\circ}{f}+\sum_{k=1}^\infty((p_k(\square)\stackrel{\circ}{f}-
p_{k-1}(\square)\stackrel{\circ}{f})t^k,$$
with $\square$ being a q-analogue of the Laplace-Beltrami operator
$$\square f \stackrel{\rm def}{=}(1-zz^*)^2 \frac{\partial^{(l)}}{\partial
z^*}\frac{\partial^{(l)}f}{\partial z}=q^2 \frac{\partial^{(r)}}{\partial
z^*}\frac{\partial^{(r)}f}{\partial z}(1-zz^*)^2,$$
with $f \in{\rm Pol}({\mathbb C})_q$ and $p_k$, $k \in{\mathbb Z}_+$, being
polynomials given by (\ref{pk}).
\end{proposition}

\medskip

 It follows from the definition of the bilinear maps $m_t$, $t \in(0,1)$,
that for all $i,j,k,l \in{\mathbb Z}_+$, $f_1,f_2 \in{\rm Pol}({\mathbb
C})_q$,
$$m_t(z^if_1,f_2)=z^im_t(f_1,f_2),$$
$$m_t(f_1,f_2z^{*l})=m_t(f_1,f_2)z^{*l},$$
$$m_t(z^{*j},z^k)=B_{q,t}(z^{*j}z^k).$$
Hence for all $i,j,k,l \in{\mathbb Z}_+$ one has
\begin{equation}\label{mtBqt}m_t((z^iz^{*j}),(z^kz^{*l}))=
z^iB_{q,t}(z^{*j}z^k)z^{*l}.
\end{equation}
We are about to deduce theorem \ref{ae} from (\ref{mtBqt}) and proposition
\ref{Bqtae}. In fact, one can easily demonstrate as in \cite[proposition
8.3]{SSV5} that
$$\square(f_2(z^*)\cdot f_1(z))=q^2 \frac{\partial^{(r)}f_2(z^*)}{\partial
z^*}(1-zz^*)^2 \frac{\partial^{(l)}f_1(z)}{\partial z}=$$
$$=\frac{\partial^{(r)}f_2(z^*)}{\partial
z^*}q^{-2}(1-(1+q^{-2})z^*z+
q^{-2}z^{*2}z^2)\frac{\partial^{(l)}f_1(z)}{\partial z}$$
for arbitrary polynomials $f_1(z)$, $f_2(z^*)$. What remains is to compare
this expression for $\square$ with the definition of $\widetilde{\square}$
and apply the fact that $\{z^iz^{*j}\}_{i,j \in{\mathbb Z}_+}$ constitute a
basis in the vector space ${\rm Pol}({\mathbb C})_q$.

\bigskip

\section{A formal associativity}

\begin{proposition}\label{m}The multiplication in ${\rm Pol}({\mathbb
C})_q[[t]]$ given by the bilinear map
$$m:{\rm Pol}({\mathbb C})_q[[t]]\times{\rm Pol}({\mathbb C})_q[[t]]\to{\rm
Pol}({\mathbb C})_q[[t]],$$
\begin{equation}\label{mf}m:\sum_{j=0}^\infty a_jt^j \times
\sum_{k=0}^\infty b_kt^k \mapsto \sum_{i=0}^\infty \left(\sum_{j+k=i}a_j*b_k
\right)t^i\;\footnotemark,
\end{equation}
\footnotetext{The outward sum clearly converges in the topological vector
space ${\rm Pol}({\mathbb C})_q[[t]]$.} with $\{a_j \}_{j \in{\mathbb
Z}_+},\{b_k \}_{k \in{\mathbb Z}_+}\in{\rm Pol}({\mathbb C})_q$, is
associative.
\end{proposition}

\medskip

 {\bf Proof.} Introduce the algebra ${\rm End}_{\mathbb C}({\mathbb
C}[z]_q)$ of all linear operators in the vector space ${\mathbb C}[z]_q$,
and the algebra ${\rm End}_{\mathbb C}({\mathbb C}[z]_q)[[t]]$ of formal
series with coefficients in ${\rm End}_{\mathbb C}({\mathbb C}[z]_q)$. To
prove our statement, it suffices to establish an isomorphism of the algebra
${\rm Pol}({\mathbb C})_q[[t]]$ equipped with the multiplication $m$ and a
subalgebra of ${\rm End}_{\mathbb C}({\mathbb C}[z]_q)[[t]]$ given the
standard multiplication. Let ${\cal I}:{\rm Pol}({\mathbb C})_q \to{\rm
End}_{\mathbb C}({\mathbb C}[z]_q)[[t]]$ be such a linear operator that for
all $j,k,m \in{\mathbb Z}_+$
$${\cal I}(z^jz^{*k}):z^m \mapsto \left
\{\begin{array}{ccl}\dfrac{(q^{2m};q^{-2})_k}{(tq^{2m};q^{-2})_k}z^{m-k+j}&,&
k \le m \\ 0 &,& k>m \end{array}\right..$$
(More precisely, one should replace the rational function
$1/(tq^{2m};q^{-2})_k$ of an indeterminate $t$ with its Teylor expansion.)

 The following lemma follows from the construction of \cite[section
7]{SSV5}.

\medskip

\begin{lemma}\label{Q}The linear map
$$Q:{\rm Pol}({\mathbb C})_q[[t]]\to{\rm End}_{\mathbb C}({\mathbb
C}[z]_q)[[t]],$$
$$Q:\sum_{j=0}^\infty f_jt^j \mapsto \sum_{j=0}^\infty{\cal
I}(f_j)t^j\;\footnote{The convergence of the series
$\sum \limits_{j=0}^\infty{\cal I}(f_j)t^j$ in the space ${\rm End}_{\mathbb
C}({\mathbb C}[z]_q)[[t]]$ is obvious.},\qquad \{f_j \}_{j \in{\mathbb
Z}_+}\subset{\rm Pol}({\mathbb C})_q,$$ is injective, and for all
$\psi_1,\psi_2 \in{\rm Pol}({\mathbb C})_q[[t]]$ one has
$Qm(\psi_1,\psi_2)=(Q \psi_1)\cdot(Q \psi_2)$.
\end{lemma}

\medskip

 Lemma \ref{Q} implies the associativity of the multiplication $m$ in
${\rm Pol}({\mathbb C})_q[[t]]$. Thus, proposition \ref{m} is proved. \hfill
$\blacksquare$

\medskip

 Define a linear operator $*$ in ${\rm Pol}({\mathbb C})_q[[t]]$ by
$$\left(\sum_{j=0}^\infty f_jt^j \right)^*=\sum_{j=0}^\infty f_j^*t^j,\qquad
\{f_j \}_{j \in{\mathbb Z}_+}\subset{\rm Pol}({\mathbb C})_q.$$

\medskip

\begin{proposition} $*$ is an involution in ${\rm Pol}({\mathbb C})_q[[t]]$
equipped by $m$ as a multiplication:
$$m(\psi_1,\psi_2)^*=m(\psi_2^*,\psi_1^*),\qquad \psi_1,\psi_2 \in{\rm
Pol}({\mathbb C})_q[[t]].$$
\end{proposition}

\smallskip

 {\bf Proof.} For all $f_1,f_2 \in{\rm Pol}({\mathbb C})_q$ one has
$$(m_0(f_1 \otimes f_2))^*=m_0(f_2^*\otimes f_1^*),$$
$$\widetilde{\square}^{21}(f_1 \otimes
f_2)^{*\otimes*}=\widetilde{\square}(f_1^*\otimes f_2^*)$$
with $\widetilde{\square}^{21}=c_0 \square c_0$, and $c_0$ being the flip of
tensor multiples. What remains is to observe that the coefficients of
$p_n(x)$, $n \in{\mathbb Z}_+$, are real, and to apply (\ref{mf}),
(\ref{stpr}), (\ref{Ck}). \hfill $\blacksquare$

\bigskip

\section{\boldmath $U_q \mathfrak{su}_{1,1}$-invariance}

 Remind some well known results on the quantum group $SU_{1,1}$ and the
quantum disc (see, for example, \cite{CP, SSV2}).

 The quantum universal enveloping algebra $U_q \mathfrak{sl}_2$ is a Hopf
algebra over ${\mathbb C}$ determined by the generators $K$, $K^{-1}$, $E$,
$F$, and the relations
$$KK^{-1}=K^{-1}K=1,\qquad K^{\pm 1}E=q^{\pm 2}EK^{\pm 1},\qquad K^{\pm
1}F=q^{\mp 2}FK^{\pm 1},$$
$$EF-FE=(K-K^{-1})/(q-q^{-1}).$$
Comultiplication $\Delta:U_q \mathfrak{sl}_2 \to U_q \mathfrak{sl}_2 \otimes
U_q \mathfrak{sl}_2$, counit $\varepsilon:U_q \mathfrak{sl}_2 \to{\mathbb
C}$ and antipode $S:U_q \mathfrak{sl}_2 \to U_q \mathfrak{sl}_2$ are given
by
$$\Delta(K^{\pm 1})=K^{\pm 1}\otimes K^{\pm 1},\qquad \Delta(E)=E \otimes
1+K \otimes E,\qquad \Delta(F)=F \otimes K^{-1}+1 \otimes F,$$
$$\varepsilon(E)=\varepsilon(F)=\varepsilon(K^{\pm 1}-1)=0,$$
$$S(K^{\pm 1})=K^{\mp 1},\qquad S(E)=-K^{-1}E,\qquad S(F)=-FK.$$
The structure of Hopf algebra allows one to define a tensor product of $U_q
\mathfrak{sl}_2$-modules and a tensor product of their morphisms. Thus, we
obtain a tensor category of $U_q \mathfrak{sl}_2$-modules.

 Consider an algebra $F$ equipped also with a structure of $U_q
\mathfrak{sl}_2$-module. $F$ is called a $U_q \mathfrak{sl}_2$-module
algebra if the multiplication
$$m_F:F \otimes F \to F,\qquad m_F:f_1 \otimes f_2 \mapsto f_1f_2,\qquad
f_1,f_2 \in F,$$
is a morphism of $U_q \mathfrak{sl}_2$-modules. (In the case $F$ has a unit,
the above definition should also include its invariance: $\xi \cdot
1=\varepsilon(\xi)1$, $\xi \in U_q \mathfrak{sl}_2$).

 The following relations determine a structure of $U_q
\mathfrak{sl}_2$-module algebra on ${\mathbb C}[z]_q$:
\begin{equation}\label{Czq}K^{\pm 1}z=q^{\pm 2}z,\qquad Fz=q^{1/2},\qquad
Ez=-q^{1/2}z^2.
\end{equation}
Equip $U_q \mathfrak{sl}_2$ with an involution:
$$E^*=-KF,\qquad F^*=-EK^{-1},\qquad (K^{\pm 1})^*=K^{\pm 1},$$
and let $U_q \mathfrak{su}_{1,1}$ stand for the Hopf $*$-algebra produced
this way. An involutive algebra $F$ is said to be $U_q
\mathfrak{su}_{1,1}$-module algebra if it is $U_q \mathfrak{sl}_2$-module
algebra, and the involutions in $F$ and $U_q \mathfrak{su}_{1,1}$ agree as
follows:
\begin{equation}\label{inv}(\xi f)^*=(S(\xi))^*f^*,\qquad \xi \in U_q
\mathfrak{su}_{1,1},\;f \in F.
\end{equation}
(\ref{Czq}) determines a structure of $U_q \mathfrak{su}_{1,1}$-module
algebra in ${\rm Pol}({\mathbb C})_q$. Thus, each of the vector spaces ${\rm
Pol}({\mathbb C})_q$, ${\rm Pol}({\mathbb C})_q[[t]]$ is equipped with a
structure of $U_q \mathfrak{su}_{1,1}$-module.

\medskip

\begin{proposition} ${\rm Pol}({\mathbb C})_q[[t]]$ with the multiplication
defined above and the involution $*$ is a $U_q \mathfrak{su}_{1,1}$-module
algebra.
\end{proposition}

\smallskip

 {\bf Proof.} Since ${\rm Pol}({\mathbb C})_q$ is a $U_q
\mathfrak{su}_{1,1}$-module algebra, (\ref{inv}) is valid for $F={\rm
Pol}({\mathbb C})_q$. Hence it is also true for $F={\rm Pol}({\mathbb
C})_q[[t]]$. What remains is to prove that ${\rm Pol}({\mathbb C})_q[[t]]$
is a $U_q \mathfrak{sl}_2$-module algebra. For that, by a virtue of
(\ref{stpr}), (\ref{Ck}), it suffices to demonstrate that the linear maps
$m_0$ and $\widetilde{\square}$ are morphisms of $U_q
\mathfrak{sl}_2$-modules. As for $m_0$, this property has already been
mentioned. So we need only to consider $\widetilde{\square}$. Given any
polynomials $f_1(z^*)$, $f_2(z)$, it follows from
$\square(f_1(z^*)f_2(z))=\sum \limits_{jk}b_{jk}z^{*j}z^k$,
$b_{jk}\in{\mathbb C}$, that $\widetilde{\square}(f_1(z^*)\otimes
f_2(z))=\sum \limits_{jk}b_{jk}z^{*j}\otimes z^k$, and
$${\widetilde{\square}}(g_1(z)f_1(z^*)\otimes f_2(z)g_2(z^*))=(g_1(z)\otimes
1)\widetilde{\square}(f_1(z^*)\otimes f_2(z))(1 \otimes g_2(z^*)).$$
Thus, it suffices to prove that $\square$ is a morphism $U_q
\mathfrak{sl}_2$-modules. This latter result is obtained in \cite{SSV2} (It
is a consequence of $U_q \mathfrak{su}_{1,1}$-invariance of the differential
calculus in the quantum disc considered there). \hfill $\blacksquare$

\medskip

{\sc Remark 6.2.} The works \cite{SSV1, SSV2} deal with the $U_q
\mathfrak{su}_{1,1}$-module algebra $D({\mathbb U})_q$ of 'finite functions
in the quantum disc'. (The space $D({\mathbb U})_q'$ mentioned in this work
is dual to $D({\mathbb U})_q$). The relations (\ref{stpr}) -- (\ref{pk})
determine a formal deformation of $D({\mathbb U})_q$ in the class of $U_q
\mathfrak{su}_{1,1}$-module algebras, that is, it allows one to equip
$D({\mathbb U})_q[[t]]$ with a structure of $U_q \mathfrak{su}_{1,1}$-module
algebra over the ring ${\mathbb C}[[t]]$.

\bigskip

\section{Concluding notes}

 We have demonstrated that the method of Berezin allows one to produce a
formal deformation for a q-analogue of the unit disc. In \cite{SV},
q-analogues for arbitrary bounded symmetric domains were constructed. We
hope in that essentially more general setting, the method of Berezin will
help remarkable results to be obtained.

\end{document}